\documentclass[12pt,reqno]{amsart}

\usepackage{verbatim}
\usepackage{color}
\usepackage{amscd}
\usepackage{amsfonts}
\usepackage{epsfig}
\usepackage{amssymb}
\usepackage{amsmath}
\usepackage{amsthm}
\usepackage{mathrsfs}
\usepackage{multicol}

\begin{document}

\numberwithin{equation}{section}
\renewcommand{\theequation}{\thesection.\arabic{equation}}
\setcounter{secnumdepth}{2}

\newcommand{\Ac}{\mathcal A}
\newcommand{\Bc}{\mathcal B}
\newcommand{\Ec}{\mathcal E}
\newcommand{\Ect}{\tilde{\Ec}}
\newcommand{\eps}{\epsilon}
\newcommand{\Fc}{\mathcal F}
\newcommand{\Gc}{\mathcal G}
\newcommand{\Hc}{\mathcal H}
\newcommand{\Ic}{\mathcal I}
\newcommand{\Lc}{\mathcal{L}_A}
\newcommand{\Tc}{\mathcal T}
\newcommand{\Mcl}{{\mathcal M}}
\newcommand{\Nb}{\mbox{$\mathbb N$}}
\newcommand{\NxN}{N \times N}
\newcommand{\Om}{\Omega}
\newcommand{\Omb}{\overline{\Omega}}
\newcommand{\pal}{\partial}
\newcommand{\pT}{p_T}
\newcommand{\Pc}{\mathcal P}
\newcommand{\Pco}{\Pc_1}
\newcommand{\Pck}{\Pc_ k}
\newcommand{\R}{\mathbb R}
\newcommand{\Rb}{\bar{\R}}
\newcommand{\Ro}{\overline{\R}}
\newcommand{\Rn}{{\R}^n}
\newcommand{\RN}{{\R}^N}
\newcommand{\Rpp}{[0,\infty)}
\newcommand{\ra}{\rightarrow}
\newcommand{\Rz}{R_0}
\newcommand{\Sc}{\mathcal{S}_A}
\newcommand{\Sct}{\tilde{\Sc}}
\newcommand{\ub}{\bar{u}}
\newcommand{\Ub}{\bar{U}}
\newcommand{\uh}{\hat{u}}
\newcommand{\UR}{U_R}
\newcommand{\ut}{\tilde{u}}
\newcommand{\vt}{\tilde{v}}
\newcommand{\vap}{\varphi}
\newcommand{\Wp}{W^{1,\hspace{0.2mm}p}}
\newcommand{\Dp}{D^{1,\hspace{0.2mm}p}}

\newcommand{\deqs}{:=}
\newcommand{\eqs}{\ =\ }
\newcommand{\geqs}{\ \geq \ }
\newcommand{\leqs}{\ \leq \ }
\newcommand{\lts}{\ < \ }
\newcommand{\n}[1]{\left\vert#1\right\vert}
\newcommand{\nm}[1]{\left\Vert#1\right\Vert}
\newcommand{\pls}{\, + \, }
\newcommand{\plms}{\ + \ }

\newcommand{\aip}[1]{\left\langle#1\right\rangle}
\newcommand{\bdU}{\partial U}
\newcommand{\bdUR}{\partial U_R}
\newcommand{\bdy}{\partial\Omega}
\newcommand{\betmu}{\beta(\mu)}
\newcommand{\beto}{\beta_1}
\newcommand{\Boc}{\overline{B_1}}
\newcommand{\delj}{\delta_j}
\newcommand{\delk}{\delta_k}
\newcommand{\deljmu}{\delta_j(\mu)}
\newcommand{\delo}{\delta_1}
\newcommand{\DFcv}{D\Fc(v)}
\newcommand{\DIk}{\pal \, I_k}
\newcommand{\Div}{\mathop{\rm div}\nolimits}
\newcommand{\Dnu}{\frac{\partial u}{\partial\nu}}
\newcommand{\dx}{\hspace{0.36mm}dx}
\newcommand{\ds}{\hspace{0.36mm}ds}
\newcommand{\dt}{\hspace{0.36mm}dt}
\newcommand{\dsg}{\hspace{0.36mm}d\sigma}
\newcommand{\acU}{\mathcal{A}\&c,\hspace{0.2mm}U}
\newcommand{\gams}{\gamma s}
\newcommand{\gau}{\gamma u}
\newcommand{\gamu}{\gamma (u)}
\newcommand{\gamst}{\gamma^*}
\newcommand{\gav}{\gamma v}
\newcommand{\grads}{\nabla s}
\newcommand{\gradu}{\nabla u}
\newcommand{\gradv}{\nabla v}
\newcommand{\gradvap}{\nabla\vap}
\newcommand{\gradw}{\nabla w}
\newcommand{\gas}{\gamma s}
\newcommand{\lamj}{\lambda_j}
\newcommand{\Iby}{\int_{\bdy}}
\newcommand{\IOm}{\int_{\Om}}
\newcommand{\InU}{\int_U}
\newcommand{\InbU}{\int_{\bdU}}
\newcommand{\Lco}{{\mathcal L}_1}
\newcommand{\Lcc}{{\mathcal L}_c}
\newcommand{\mns}{\, - \, }
\newcommand{\nei}{\not\in}
\newcommand{\nmc}{\left|\left|\cdot\right|\right|}
\newcommand{\nmu}{\left|\left|u\right|\right|}
\newcommand{\Rp}{[0,\infty)}
\newcommand{\sot}{_{1,2}}
\newcommand{\sjt}{ \tilde{s}_j }
\newcommand{\skt}{ \tilde{s}_k }
\newcommand{\uk}{ \mathfrak{u}_k }
\newcommand{\vk}{ \mathfrak{v}_k }

\newcommand{\foral}{\quad \mbox{for all} \quad}
\newcommand{\form}{\quad \mbox{for some} \quad}
\newcommand{\Stek}{\mbox{Steklov eigenfunction of} \ \Lco \ \mbox{on} \ U}
\newcommand{\StekL}{\mbox{Steklov eigenfunction of} \ \left(\Lc, b\right) }
\newcommand{\xand}{\ \mbox{and} \quad }
\newcommand{\xfor}{ \quad \mbox{for} \ }
\newcommand{\xiff}{\ \mbox{if and only if} \ }
\newcommand{\xon}{\qquad \mbox{on} \quad}
\newcommand{\xor}{\qquad \mbox{or} \quad}
\newcommand{\wrt}{\mbox{with respect to}}
\newcommand{\xwlog}{\mbox{without loss of generality}}

\newcommand{\barr}{\begin{eqnarray}}
\newcommand{\beq}{\begin{equation}}
\newcommand{\bpf}{\begin{proof}\quad}
\newcommand{\btm}{\begin{thm}}
\newcommand{\bp}{\begin{prop}}
\newcommand{\earr}{\end{eqnarray}}
\newcommand{\eeq}{\end{equation}}
\newcommand{\epf}{\end{proof}}
\newcommand{\etm}{\end{thm}}
\newcommand{\ep}{\end{prop}}

\newcommand{\Cbdy}{C(\bdy)}
\newcommand{\CcOU}{C_c^1\left(U\right)}
\newcommand{\CCON}{C^1_c\left(\RN\right)}
\newcommand{\COm}{C\left(\Omb\right)}
\newcommand{\Coneb}{C^1\left(\Omb\right)}
\newcommand{\CROm}{C^1_R\left(\Omb\right)}
\newcommand{\DoN}{D^1\left(\RN\right)}
\newcommand{\DpN}{D^{1,\hspace{0.2mm}p}\left(\RN\right)}
\newcommand{\EoUk}{E^1_k\left(U\right)}
\newcommand{\EoU}{E^1\left(U\right)}
\newcommand{\EoUo}{E_0^1\left(U\right)}
\newcommand{\EoUs}{\EoU \times \EoU}
\newcommand{\Gop}{G^{1,\hspace{0.2mm}p}\left(\Om\right)}
\newcommand{\HarmU}{\mathscr{H}\left(U\right)}
\newcommand{\Hhe}{H^{1/2}_e\left(\bdU\right)}
\newcommand{\Hhed}{H^{-1/2}_e\left(\bdU\right)}
\newcommand{\Hhb}{H^{1/2}\left(\bdy\right)}
\newcommand{\HhU}{H^{1/2}\left(\bdU\right)}
\newcommand{\HhUd}{H^{-1/2}\left(\bdU\right)}
\newcommand{\Hsb}{H^s\left(\bdy\right)}
\newcommand{\HsbUe}{H^s_e\left(\bdU\right)}
\newcommand{\HR}{H^1\left(\RN\right)}
\newcommand{\Hone}{H^1\left(\Om\right)}
\newcommand{\Ho}{H^1_{\widehat{\mathrm{o}}}}
\newcommand{\HoU}{H^1\left(U\right)}
\newcommand{\HoUs}{\HoU \times \HoU}
\newcommand{\HoUo}{H_0^1\left(U\right)}
\newcommand{\HoG}{H^1\left(G\right)}
\newcommand{\HoGo}{H_0^1\left(G\right)}
\newcommand{\Hones}{\Hone \times \Hone}
\newcommand{\Hz}{H_{0}^1\left(\Om\right)}
\newcommand{\HarU}{\mathcal{H}_A\left(U\right)}
\newcommand{\Lb}{L^p\left(\pal\Om,\dsg\right)}
\newcommand{\Linb}{L^{\infty}\left(\bdy,\dsg\right)}
\newcommand{\LinbU}{L^{\infty}\left(\bdU,\dsg\right)}
\newcommand{\Lpb}{L^p\left(\bdy,\dsg\right)}
\newcommand{\LplU}{L^p_{loc}\left(U\right)}
\newcommand{\LpTb}{L^{p_T}\left(\bdy,\dsg\right)}
\newcommand{\LpbU}{L^p\left(\bdU,\dsg\right)}
\newcommand{\Lqb}{L^q\left(\bdy,\dsg\right)}
\newcommand{\LqbU}{L^q\left(\bdU,\dsg\right)}
\newcommand{\LqTb}{L^{q_T}\left(\bdy,\dsg\right)}
\newcommand{\LrbU}{L^r\left(\bdU,\dsg\right)}
\newcommand{\LsbU}{L^s\left(\bdU,\dsg\right)}
\newcommand{\Ltbg}{L^2 (\bdy,\, g \, \dsg)}
\newcommand{\LtbU}{L^2\left(\bdU,\dsg\right)}
\newcommand{\LtbbyU}{L^2\left(\bdU, b \, \dsg \right)}
\newcommand{\Lp}{L^p\left(\Om\right)}
\newcommand{\LpU}{L^p\left(U\right)}
\newcommand{\Lq}{L^q\left(\Om\right)}
\newcommand{\LqU}{L^q\left(U\right)}
\newcommand{\Lt}{L^2\left(\Om\right)}
\newcommand{\LtU}{L^2\left(U\right)}
\newcommand{\NL}{N\left(\Lc\right)}
\newcommand{\NLc}{N\left(\Lcc\right)}
\newcommand{\No}{N\left(\mathcal{L}_1\right)}
\newcommand{\N}{N\left(\mathcal{L}_1\right)}
\newcommand{\NrhU}{\mathscr{N}\left(U\right)}
\newcommand{\RxHone}{\R\times\Hone}
\newcommand{\Wone}{W^{1,\hspace{0.2mm}1}\left(\Om\right)}
\newcommand{\Wop}{W^{1,\hspace{0.2mm}p}\left(\Om\right)}
\newcommand{\WpU}{W^{1,\hspace{0.2mm}p}\left(U\right)}
\newcommand{\EpU}{E^{1,\hspace{0.2mm}p}\left(U\right)}
\newcommand{\EopU}{E^{1,\hspace{0.2mm}p}_0\left(U\right)}

\newtheorem{thm}{Theorem}[section]
\newtheorem{corr}[thm]{Corollary}
\newtheorem{cond}{Condition}
\newtheorem{lem}[thm]{Lemma}
\newtheorem{prop}[thm]{Proposition}
\newtheorem{remak}[thm]{Remark}

\title[\sl Laplace's equation on exterior regions]{\sl Laplace's equation with concave and convex boundary nonlinearities on an exterior region}

\author[\tt Jinxiu Mao]{\tt Jinxiu Mao}

\address{School of Mathematics, Qufu Normal University, Qufu, Shandong 273165, China\vskip 2pt{\sf Email: maojinxiu1982@163.com}}

\author{Zengqin Zhao}
\address{School of Mathematics,
 Qufu Normal University, Qufu,
 Shandong, 273165, 
China}
\email{zqzhaoy@163.com}
\thanks{Supported by the Natural Science Foundation of China(11571197), Natural Science Foundation of Shandong Province(ZR2014AM007) and the Science
Foundation of Qufu Normal University of China (XJ201112).}

\begin{abstract}
  This paper studies Laplace's equation $-\Delta\,u=0$ in an exterior region $U\varsubsetneq\RN$, when $N\geq3$, subject to the nonlinear boundary condition $\Dnu=\lambda\n{u}^{q-2}u+\mu\n{u}^{p-2}u$ on $\bdU$ with $1<q<2<p<2^*$.
  In the function space $\HarmU$, one observes when $\lambda>0$ and $\mu\in\R$ arbitrary, then there exists a sequence $\left\{u_k\right\}$ of solutions with negative energy converging to $0$ as $k\to\infty$; on the other hand, when $\lambda\in\R$ and $\mu>0$ arbitrary, then there exists a sequence $\left\{\ut_k\right\}$ of solutions with positive and unbounded energy.
  Also, associated with the $p$-Laplacian equation $-\Delta_p\,u=0$, the exterior $p$-harmonic Steklov eigenvalue problems are described.
\end{abstract}

\maketitle

\section{Introduction}\label{Introduction}
This paper discusses the existence of infinitely many harmonic functions in an exterior region $U\varsubsetneq\RN$ when $N\geq3$, subject to a nonlinear boundary condition on $\bdU$ that combines concave and convex terms with $1<q<2<p<2^*$, described as below
\begin{equation}\label{e1.1}
\left\{\begin{array}{ll}
-\Delta\,u(x)\,=\,0&\mathrm{in}\hspace{2mm}U,\\ \\
\Dnu(z)\,=\,\lambda\n{u(z)}^{q-2}u(z)\,+\,\mu\n{u(z)}^{p-2}u(z)&\mathrm{on}\hspace{2mm}\bdU,
\end{array}\right.
\end{equation}
in the space $\EoU$ of functions where $u\in L^{2^*}\left(U\right)$ and $\gradu\in L^2\left(U;\RN\right)$.
Here, $2^*:=\frac{2N}{N-2}$ is the critical Sobolev index and $\gradu:=\left(D_1u,D_2u,\ldots,D_Nu\right)$ is the weak gradient of $u$.

A \textsl{region} is a nonempty, open, connected subset $U$ of $\RN$, and is said to be an \textsl{exterior region} provided that its complement $\RN\setminus U$ is a nonempty, compact subset.
Without loss of generality, we simply assume that $0\notin U$.
The boundary of a set $A$ is denoted by $\pal A$.

Our general assumption on $U$ is the following condition.

\vskip 4pt
\noindent{\bf Condition B1.} \textit{$U\varsubsetneq\RN$ is an exterior region, with $0\notin U$, whose boundary $\bdU$ is the union of finitely many disjoint, closed, Lipschitz surfaces, each of finite surface area.}
\vskip 4pt

One may want to notice here that the following prototypical problem
\begin{equation}
-\Delta\,u(x)\,=\,\lambda\n{u(x)}^{q-2}u(x)\,+\,\mu\n{u(x)}^{p-2}u(x)\hspace{4mm}\mathrm{in}\hspace{2mm}\Om\nonumber
\end{equation}
has originally been investigated by Ambrosetti, Br\'{e}zis and Cerami \cite{ABC} in 1994 and then in 1995 by Bartsch and Willem \cite{BW}, in the function space $\Hz$ on a bounded region $\Om$ with a smooth boundary $\bdy$.
Since then, there has been a large number of papers appearing on some related problems; nevertheless, it seems to me the description of the existence of solutions to problem \eqref{e1.1} is missing as a reasonable decomposition result of the associated Hilbert function space is required for application of the dual fountain theorem, as discussed in \cite{BW, Wi}.

The aim of this paper is to solve problem \eqref{e1.1} using a recent decomposition result by Auchmuty and Han \cite{AH2, Ha1}.
To state our result, we first define the energy functional
\begin{equation}\label{e1.2}
\varphi(u)\,:=\,\frac{1}{2}\InU\n{\gradu}^2\dx\,-\,\frac{\lambda}{q}\InbU\n{u}^q\dsg\,-\,\frac{\mu}{p}\InbU\n{u}^p\dsg.
\end{equation}
Here, $dx$ is the Lebesgue volume element of $\RN$ while $d\sigma$ is the Hausdorff $\left(N-1\right)$-dimensional surface element of $\bdU$.
The main result of this paper is described as below.

\begin{thm}\label{T1.1}
Assume condition {\bf(B1)} holds and $1<q<2<p<2^*$.\\
{\bf(a)} When $\lambda\in\R$ and $\mu>0$ arbitrary, then problem \eqref{e1.1} has a sequence $\left\{\uk\right\}$ of solutions in $\EoU$ such that $\varphi(\uk)>0$ and $\varphi(\uk)\to\infty$ as $k\to\infty$.\\
{\bf(b)} When $\lambda>0$ and $\mu\in\R$ arbitrary, then problem \eqref{e1.1} has a sequence $\left\{\vk\right\}$ of solutions in $\EoU$ such that $\varphi(\vk)<0$ and $\varphi(\vk)\to0^-$ as $k\to\infty$.
\end{thm}

We recall in section 2 some necessary results to carry out the proofs that are detailed in section 3; the last section 4 is devoted to the description of the $p$-harmonic Steklov eigenvalue problems on an exterior region $U$ in a Banach space $\EpU$ when $1<p<N$.

It is interesting to see that some nice properties of the first exterior $p$-harmonic Steklov eigenvalue problem are described in Han \cite{Ha1,Ha2,Ha4} and he \cite{Ha3} also studied an exterior harmonic boundary value problem with some oscillating boundary condition.
However, as far as I know, there is no result for the sequence of $p$-harmonic Steklov eigenvalue problems on an exterior region $U$, so I will study this in section 4.
See Torn\'{e} \cite{To} for bounded region case.

Finally, one notices that the author only wants to simply present an application of some result in \cite{AH2, Ha1}.
Theorem 1.1 remains true when the special nonlinearity in \eqref{e1.1} is replaced by more general ones as mentioned in \cite{BW}.
On the other hand, it is very interesting to know more results like this using the fountain theorems of Yan and Yang \cite{YY}, and Zou \cite{Zo}.

\vspace{2ex}
\section{The function space $\EpU$}\label{FunctionSpaces}
First, let's fix the notations that will be used in this paper.
Given $p,q\in\left[1,\infty\right]$, $\LpU$ and $\LqbU$ are the usual spaces of extended, real-valued, Lebesgue measurable functions on $U$ and $\bdU$, with their standard norms written as $\nm{\cdot}_{p,\hspace{0.2mm}U}$ and $\nm{\cdot}_{q,\hspace{0.2mm}\bdU}$, respectively.

Auchmuty and Han \cite{AH2,AH3,Ha1} recently introduced a new function space $\EpU$ suitale for the study of harmonic boundary value problems on an exterior region $U$ which satisfies the boundary regularity condition {\bf(B1)} - that is, each function $u\in\EpU$ satisfies $u\in L^{p^*}\left(U\right)$ and $\n{\gradu}\in L^p\left(U\right)$ with $N\geq2$ and $p^*:=\frac{Np}{N-p}$ when $1<p<N$.

The gradient $L^p$-norm provides a norm to guarantee $\EpU$ a Banach function space - that is, $\EpU$ is a Banach function space with respect to the norm
\begin{equation}\label{e2.1}
\nm{u}_{p,\hspace{0.2mm}\nabla}\,:=\,\left(\InU\n{\gradu}^p\dx\right)^{\frac{1}{p}}\hspace{7mm}\mathrm{for}\hspace{2mm}\mathrm{all}\hspace{2mm}u\in\EpU.
\end{equation}

Notice when $p\geq N$, Auchmuty and Han \cite{AH3} showed, with an interesting example, that $\EpU$ is not complete with respect to the gradient $L^p$-norm in general.

When $N\geq3$ and $p=2$, they instead used the notation $\EoU$ to denote the associated Hilbert function space with respect to the gradient $L^2$-inner product
\begin{equation}\label{e2.2}
\aip{u,v}_{\nabla}\,:=\,\InU\,\gradu\cdot\gradv\,\dx\hspace{7mm}\mathrm{for}\hspace{2mm}\mathrm{all}\hspace{2mm}u,\,v\in\EoU,
\end{equation}
whose norm is thus written as $\nm{u}_{\nabla}$.
In addition, one has the direct sum
\begin{equation}\label{e2.3}
\EoU\,=\,\EoUo\,\oplus\,\HarmU,
\end{equation}
where $\HarmU$ denotes the Hilbert subspace of $\EoU$ of all functions $u$ satisfying
\begin{equation}
\aip{u,v}_{\nabla}\,=\,\InU\,\gradu\cdot\gradv\,\dx\,=\,0\hspace{7mm}\mathrm{for}\hspace{2mm}\mathrm{all}\hspace{2mm}v\in\CcOU\nonumber
\end{equation}
and $\EoUo$ denotes the closure of $\CcOU$ with respect to this $\nabla$-norm.
Here, $\CcOU$ is the set of functions that are continuously differentiable and have compact support in $U$.

Let's recall some results about the space $\EpU$ which will be used later.

\vskip 4pt
\noindent\textbf{Result 2.1.} \textit{Suppose that $N\geq2$, $1<p<N$ and condition {\bf(B1)} holds.
Then, the embedding of $\EpU$ into $L^{p^*}\left(U\right)$ is continuous, where $p^*:=\frac{Np}{N-p}$ is the critical Sobolev index; besides, the embedding of $\EpU$ into $\LqbU$ is continuous when $1\leq q\leq p_*$ and also compact when $1\leq q<p_*$, where $p_*:=\frac{\left(N-1\right)p}{N-p}$ is the trace critical Sobolev index.}
\vskip 4pt

Obviously, result 2.1 shows us some concrete function spaces that are contained in the dual space of $\EpU$.
The preceding results can be found, with details, in \cite{AH2, AH3, Ha1}.

Below, we give the fountain theorems.
Given a compact group $\mathfrak{G}$ and a normed vector space $\mathcal{X}$ with norm $\nm{\cdot}$, we say $\mathfrak{G}$ acts \textsl{isometrically} on $\mathcal{X}$ provided $\nm{gu}=\nm{u}$ for all $g\in\mathfrak{G}$ and $u\in\mathcal{X}$; also, a subset $\tilde{\mathcal{X}}\subseteq\mathcal{X}$ is said to be \textsl{invariant} with respect to $\mathfrak{G}$ provided $gu\in\tilde{\mathcal{X}}$ for every $u\in\tilde{\mathcal{X}}$ and $g\in\mathfrak{G}$.
On the other hand, given $\mathfrak{G}$ and a finite dimensional space $\mathbf{V}$, we say the action of $\mathfrak{G}$ on $\mathbf{V}$ is \textsl{admissible} when each continuous, equivariant map $\wp:\pal\mathbf{O}\to\mathbf{V}^{k}$ has a zero, where $\mathbf{O}$ is an open, bounded, invariant (with respect to $\mathfrak{G}$) neighborhood of $0$ in $\mathbf{V}^{k+1}$ for some $k\geq1$; here, the map $\wp$ is said to be \textsl{equivariant} provided $g\circ\wp=\wp\circ g$ for all $g\in\mathfrak{G}$, with $g\left(v_1,v_2,\ldots,v_k\right):=\left(gv_1,gv_2,\ldots,gv_k\right)$ for any $v=\left(v_1,v_2,\ldots,v_k\right)\in\mathbf{V}^{k}$.

Next, given a Banach space $\mathcal{X}$, a functional $\psi:\mathcal{X}\to\R$ is said to belong to $C^1\left(\mathcal{X},\R\right)$, provided its first \textsl{Fr\'{e}chet derivative} exists and is continuous on $\mathcal{X}$; when $\psi$ has a continuous first \textsl{Gateaux derivative} $\psi^{\prime}$ on $\mathcal{X}$, then one observes $\psi\in C^1\left(\mathcal{X},\R\right)$.
Clearly, the functional $\varphi$ defined in \eqref{e1.2} is in $C^1\left(\HarmU,\R\right)$ and we shall assume this from now on.
Also, $\psi:\mathcal{X}\to\R$ is said to be \textsl{invariant} with respect to $\mathfrak{G}$ provided $\psi\circ g=\psi$ for every $g\in\mathfrak{G}$.

Now, let $\mathcal{X}$ be a Banach space with $\mathcal{X}=\overline{\bigoplus\limits_{j\in\mathbb{N}}\mathcal{X}(j)}$ and write, for each $k\in\mathbb{N}$,
\begin{equation}\label{e2.6}
\mathcal{Y}_k\,:=\,\overline{\oplus^k_{j=0}\,\mathcal{X}(j)}\hspace{6mm}\mathrm{and}\hspace{6mm}\mathcal{Z}_k\,:=\,\overline{\oplus^{\infty}_{j=k}\,\mathcal{X}(j)}.
\end{equation}
Then, one has the following results - \textsl{fountain theorem} and \textsl{dual fountain theorem}.

\begin{thm}\label{T2.3}
Let $\mathfrak{G}$ be a compact group, $\mathcal{X}=\overline{\bigoplus\limits_{j\in\mathbb{N}}\mathcal{X}(j)}$ a Banach space with norm $\nm{\cdot}$, and $\psi\in C^1\left(\mathcal{X},\R\right)$ an invariant functional; for each $k\in\mathbb{N}$, let $\mathcal{Y}_k,\mathcal{Z}_k$ be defined as in \eqref{e2.6}, and $\rho_k>\varrho_k>0$ some constants.
For every $k\geq k_0$ with a fixed $k_0\in\mathbb{N}$, we also assume\\
{\bf(a1)} $\mathfrak{G}$ acts isometrically on $\mathcal{X}$, the spaces $\mathcal{X}(j)$ are invariant and there is a finite dimensional space $\mathbf{V}$ such that, for all $j\in\mathbb{N}$, $\mathcal{X}(j)\simeq\mathbf{V}$ and the action of $\mathfrak{G}$ on $\mathbf{V}$ is admissible;\\
{\bf(a2)} $\mathfrak{a}_k\,:=\,\max\limits_{u\in\mathcal{Y}_k,\,\nm{u}=\rho_k}\psi(u)\,<\,0$;\\
{\bf(a3)} $\mathfrak{b}_k\,:=\,\inf\limits_{u\in\mathcal{Z}_k,\,\nm{u}=\varrho_k}\psi(u)\,\to\,\infty$ as $k\to\infty$;\\
{\bf(a4)} $\psi$ satisfies the $(PS)_c$-condition for every $c\in\left(0,\infty\right)$.\\
Then, $\psi$ has a sequence of critical values $\left\{\uk\right\}$ with $\psi(\uk)>0$ and $\psi(\uk)\to\infty$ when $k\to\infty$.
\end{thm}

\begin{thm}\label{T2.4}
Under the hypotheses of theorem \ref{T2.3}, suppose again that condition {\bf(a1)} holds.
For every $k\geq k_1$ with a fixed $k_1\in\mathbb{N}$, we also assume\\
{\bf(b1)} $\tilde{\mathfrak{a}}_k\,:=\,\max\limits_{u\in\mathcal{Y}_k,\,\nm{u}=\varrho_k}\psi(u)\,<\,0$;\\
{\bf(b2)} $\tilde{\mathfrak{b}}_k\,:=\,\inf\limits_{u\in\mathcal{Z}_k,\,\nm{u}=\rho_k}\psi(u)\,\geq\,0$;\\
{\bf(b3)} $\tilde{\mathfrak{c}}_k\,:=\,\inf\limits_{u\in\mathcal{Z}_k,\,\nm{u}\leq\rho_k}\psi(u)\,\to\,0^-$ as $k\to\infty$;\\
{\bf(b4)} $\psi$ satisfies the $(PS)^*_c$-condition with respect to $\mathcal{Y}_k$ for each $c\in\left[\tilde{\mathfrak{c}}_{k_1},0\right)$.\\
Then, $\psi$ has a sequence of critical values $\left\{\vk\right\}$ with $\psi(\vk)<0$ and $\psi(\vk)\to0^-$ when $k\to\infty$.
\end{thm}

\vskip 4pt
\noindent\textsf{Remark.} Notice $\tilde{\mathfrak{c}}_k\leq\min\limits_{u\in\mathcal{X}(k),\,\nm{u}=\varrho_k}\psi(u)\leq\max\limits_{u\in\mathcal{X}(k),\,\nm{u}=\varrho_k}\psi(u)\leq\tilde{\mathfrak{a}}_k<0$ as $\mathcal{Y}_k\cap\mathcal{Z}_k=\mathcal{X}(k)$ - this fact is used in conditions {\bf(b3)} and {\bf(b4)} presented above in theorem \ref{T2.4}.
\vskip 4pt

A sequence $\left\{u_k\right\}$ is said to be a \textsl{Palais-Smale sequence} for the functional $\psi\in C^1\left(\mathcal{X},\R\right)$ at level $c$ in $\mathcal{X}$, \textsl{$(PS)_c$-sequence} for short, if $\psi(u_k)\to c$ yet $\psi^{\prime}(u_k)\to0$ as $k\to\infty$; $\psi$ satisfies the \textsl{$(PS)_c$-condition} provided each $(PS)_c$-sequence has a strongly convergent subsequence in $\mathcal{X}$.
On the other hand, a sequence $\left\{\ut_{k_l}\right\}$, with $\ut_{k_l}$ in $\mathcal{Y}_{k_l}$, is said to be a \textsl{generalized Palais-Smale sequence} for $\psi$ at level $c$, \textsl{$(PS)^*_c$-sequence} for short,
if $\psi(\ut_{k_l})\to c$ yet $\psi|^{\prime}_{\mathcal{Y}_{k_l}}(\ut_{k_l})\to0$ as $l\to\infty$; $\psi$ satisfies the \textsl{$(PS)^*_c$-condition with respect to $\mathcal{Y}_k$} provided each $(PS)^*_c$-sequence has a subsequence that converges strongly to a critical point of $\psi$ in $\mathcal{X}$.

More details on fountain theorems can be found in \cite{Ba,BW,Wi,YY,Zo}.

\vspace{2ex}
\section{Existence results of \eqref{e1.1}}\label{Proofs}
In this section, we shall present the proofs of theorem \ref{T1.1}.
Matching with the preceding notations, we can identify $\mathfrak{G}=\mathbb{Z}_2$ - the second order quotient group, $\mathcal{X}=\HarmU$ - the Hilbert subspace of $\EoU$ of all finite energy harmonic functions, and $\psi=\varphi\in C^1\left(\mathcal{X},\R\right)$.
One result in \cite{AH2} shows $\mathcal{X}=\overline{\bigoplus\limits_{j\in\mathbb{N}}\mathcal{X}(j)}$; here, $\mathcal{X}(j)=\mathrm{span}\left\{s_j\right\}\simeq\mathbf{V}=\R$, with $s_j\in\HarmU$ a finite energy \textsl{harmonic Steklov eigenfunction} associated with the $j$-th \textsl{harmonic Steklov eigenvalue} $\delta_j>0$.
Noting that the functional $\varphi$ is even, condition {\bf(a1)} is trivially satisfied since a classical result of Borsuk-Ulam says that the antipodal action of $\mathbb{Z}_2$ on $\R$ is admissible.

In the following, we shall deduce conditions {\bf(a2)}-{\bf(a4)} and {\bf(b1)}-{\bf(b4)} to guarantee the conclusions of the first and the second part of theorem \ref{T1.1}, respectively.

To further simplify notations, set $\nm{\cdot}:=\nm{\cdot}_{\nabla}$ and $\nm{\cdot}_s:=\nm{\cdot}_{s,\hspace{0.2mm}\bdU}$ in this section.
Take a nonzero $u\in\mathcal{Y}_k$, and use $tu$ and \eqref{e1.2} to derive, for some sufficiently large $t>0$,
\begin{equation}\label{e3.1}
\varphi(tu)\,:=\,\frac{t^2}{2}\nm{u}^2\,-\,\frac{\lambda\,t^q}{q}\nm{u}^q_q\,-\,\frac{\mu\,t^p}{p}\nm{u}^p_p\,<\,0
\end{equation}
in view of the fact that the space $\mathcal{Y}_k$ is of finite dimension - so that all norms are equivalent.
As such, condition {\bf(a2)} is satisfied for every $\rho_k\geq t\nm{u}>0$ when $\mu>0$.

Next, define
\begin{equation}\label{e3.2}
\alpha_k\,:=\,\sup\limits_{u\in\mathcal{Z}_k-\{0\}}\frac{\nm{u}_p}{\nm{u}}\,>\,0.
\end{equation}
Then, one observes $\alpha_k\to0$ when $k\to\infty$.
Actually, it is readily seen that $0<\alpha_{k+1}\leq\alpha_k$, so that $\alpha_k\to\alpha\geq0$ as $k\to\infty$.
By hypotheses, there exists a $u_k\in\mathcal{Z}_k$ satisfying $\nm{u_k}=1$ and $\nm{u_k}_p\geq\alpha_k/2$ for each $k$; by definition of $\mathcal{Z}_k$, one sees $u_k\rightharpoonup0$, \textit{i.e.}, $u_k$ converges weakly to $0$, in $\HarmU$.
Result 2.1 then yields $u_{k_l}\to0$ in $\LpbU$ as $l\to\infty$, for a subsequence $\left\{u_{k_l}\right\}$ of $\left\{u_k\right\}$.
That is, $\alpha=0$.
On each subspace $\mathcal{Z}_k$ with a sufficiently large norm, we have
\begin{equation}\label{e3.3}
\varphi(u)\,\geq\,\frac{1}{2}\nm{u}^2\,-\,\frac{\n{\lambda}c^q_1}{q}\nm{u}^q\,-\,\frac{\mu\,\alpha^p_k}{p}\nm{u}^p
\,\geq\,\frac{1}{2}\left(\frac{1}{2}+\frac{1}{p}\right)\nm{u}^2\,-\,\frac{\mu\,\alpha^p_k}{p}\nm{u}^p,
\end{equation}
where $c_1>0$ is such a constant that $\nm{u}_q\leq c_1\nm{u}$ for all $u\in\EoU$.
Take $\varrho_k:=\left(\mu\alpha^p_k\right)^{-\frac{1}{p-2}}$.
Then, via the property of $\alpha_k$, there exists a $k_0\in\mathbb{N}$ such that one can always choose $\rho_k=2\varrho_k$ in the foregoing estimate \eqref{e3.1} for every $k\geq k_0$; in addition, for \eqref{e3.3}, one derives
\begin{equation}\label{e3.4}
\varphi(u)\,\geq\,\frac{1}{2}\left(\frac{1}{2}-\frac{1}{p}\right)\mu^{-\frac{2}{p-2}}\,\alpha^{-\frac{2p}{p-2}}_k\,\to\,\infty
\end{equation}
when $k\to\infty$.
As a consequence, condition {\bf(a3)} is assured.

Finally, take a $(PS)_c$-sequence $\left\{u_k\right\}$ for the functional $\varphi$ at level $c>0$ in $\HarmU$; that is, $\varphi(u_k)\to c$ yet $\varphi^{\prime}(u_k)\to0$ as $k\to\infty$.
Then, one has, for $k$ sufficiently large,
\begin{equation}\label{e3.5}
\begin{split}
c+1+\nm{u_k}&\,\geq\,\varphi(u_k)\,-\,\frac{1}{p}\,\varphi^{\prime}(u_k)(u_k)\\
&\,=\,\left(\frac{1}{2}-\frac{1}{p}\right)\InU\n{\gradu_k}^2\dx\,-\,\lambda\left(\frac{1}{q}-\frac{1}{p}\right)\InbU\n{u_k}^q\dsg,
\end{split}
\end{equation}
from which one deduces immediately that
\begin{equation}\label{e3.6}
c+1+\nm{u_k}+\n{\lambda}c^q_1\left(\frac{1}{q}-\frac{1}{p}\right)\nm{u_k}^q\,\geq\,\left(\frac{1}{2}-\frac{1}{p}\right)\nm{u_k}^2.
\end{equation}
As a result, these $u_k$'s are bounded, and thus, without loss of generality, converge weekly to a function $u\in\HarmU$; via result 2.1 again, we may simply suppose that $u_k\to u$ in $\LpbU$ and $\LqbU$ when $k\to\infty$; besides, using \eqref{e1.2}, a routine calculation leads to
\begin{equation}\label{e3.7}
\begin{split}
\nm{u_k-u}^2\,=\,&\left(\varphi^{\prime}(u_k)-\varphi^{\prime}(u)\right)(u_k-u)\,+\,\lambda\InbU\left[\,\n{u_k}^{q-2}u_k-\n{u}^{q-2}u\,\right](u_k-u)\,\dsg\\
&\,+\,\mu\InbU\left[\,\n{u_k}^{p-2}u_k-\n{u}^{p-2}u\,\right](u_k-u)\,\dsg\,\to\,0
\end{split}
\end{equation}
when $k\to\infty$.
Thus, condition {\bf(a4)} is also derived so that part {\bf(a)} of theorem \ref{T1.1} is proved.

On the other hand, apply a parallel idea as shown in \eqref{e3.1} to prove that condition {\bf(b1)} is satisfied for every $0<\varrho_k\leq t\nm{u}$, with a nonzero $u\in\mathcal{Y}_k$ and some sufficiently small $t>0$, when $\lambda>0$ - since $\mathcal{Y}_k$ is of finite dimension.
Next, define
\begin{equation}\label{e3.8}
\beta_k\,:=\,\sup\limits_{u\in\mathcal{Z}_k-\{0\}}\frac{\nm{u}_q}{\nm{u}}\,>\,0.
\end{equation}
Similarly, one observes that $\beta_k\to0$ when $k\to\infty$ again from result 2.1.
On each subspace $\mathcal{Z}_k$ with a sufficiently small norm, we have
\begin{equation}\label{e3.9}
\varphi(u)\,\geq\,\frac{1}{2}\nm{u}^2\,-\,\frac{\lambda\,\beta^q_k}{q}\nm{u}^q\,-\,\frac{\n{\mu}c^p_2}{p}\nm{u}^p\,\geq\,\frac{1}{4}\nm{u}^2\,-\,\frac{\lambda\,\beta^q_k}{q}\nm{u}^q,
\end{equation}
where $c_2>0$ is such a constant that $\nm{u}_p\leq c_2\nm{u}$ for any $u\in\EoU$.
Take $\rho_k:=\left(\frac{4\lambda\beta^q_k}{q}\right)^{\frac{1}{2-q}}$ to derive $\varphi(u)\geq0$.
Via the property of $\beta_k$, there exists a $k_1\in\mathbb{N}$ such that we may select the above $\varrho_k=\frac{\rho_k}{2}$ to be sufficiently small for all $k\geq k_1$.
As such, condition {\bf(b2)} is assured; also, condition {\bf(b3)} follows in view of the fact $\rho_k\to0$ as $k\to\infty$.
Finally, take a $(PS)^*_c$-sequence $\left\{\ut_{k_l}\right\}$, with $\ut_{k_l}$ in $\mathcal{Y}_{k_l}$, for $\varphi$ at level $c\in\left[\tilde{\mathfrak{c}}_{k_1},0\right)$; that is, $\varphi(\ut_{k_l})\to c$ whereas $\varphi|^{\prime}_{\mathcal{Y}_{k_l}}(\ut_{k_l})\to0$ as $l\to\infty$.
So, one infers $\frac{1}{p}\n{\varphi^{\prime}(\ut_{k_l})(\ut_{k_l})}\leq\nm{\ut_{k_l}}$ for sufficiently large $l$.
Thus, \eqref{e3.6} holds again.
As a result, these $\ut_{k_l}$'s are bounded, and thus converge weekly without loss of generality to a function $\ut\in\HarmU$; via result 2.1 again, we may simply assume that $\ut_{k_l}\to\ut$ in $\LpbU$ and $\LqbU$ as $l\to\infty$; by use of \eqref{e1.2} again, one observes, just like \eqref{e3.7}, $\nm{\ut_{k_l}-\ut}\to0$ when $l\to\infty$.
Thus, condition {\bf(b4)} is also derived so that part {\bf(b)} of theorem \ref{T1.1} is proved.

All the above discussions finishes the proof of theorem \ref{T1.1} completely. $\hspace{26mm}\blacksquare$

We do not know whether $\vk\to0$ as $k\to\infty$; this is the case if $0$ is the only solution of problem \eqref{e1.1} with energy $0$.
However, we can derive the following result.

\begin{prop}\label{P3.1}
Assume condition {\bf(B1)} holds and $1<q<2<p<2^*$.\\
{\bf(a)} When $\lambda\in\R$ arbitrary yet $\mu\leq0$, then \eqref{e1.1} has no solution with positive energy; also,
\begin{equation}
\inf\left\{\,\nm{u}:u\hspace{2mm}solves\hspace{2mm}\eqref{e1.1}\hspace{2mm}with\hspace{2mm}\varphi(u)>0\,\right\}\,\to\,\infty\hspace{2mm}as\hspace{2mm}\mu\,\to\,0^+.\nonumber
\end{equation}
{\bf(b)} When $\mu\in\R$ arbitrary yet $\lambda\leq0$, then \eqref{e1.1} has no solution with negative energy; also,
\begin{equation}
\sup\left\{\,\nm{v}:v\hspace{2mm}solves\hspace{2mm}\eqref{e1.1}\hspace{2mm}with\hspace{2mm}\varphi(v)<0\,\right\}\,\to\,0\hspace{2mm}as\hspace{2mm}\lambda\,\to\,0^+.\nonumber
\end{equation}
\end{prop}

\begin{proof}
Take $u\in\HarmU$ to be such that $\varphi(u)\geq0$ and $\varphi^{\prime}(u)=0$.
Then, one has
\begin{equation}\label{e3.10}
\varphi(u)\,-\,\frac{1}{q}\,\varphi^{\prime}(u)(u)\,=\,\left(\frac{1}{2}-\frac{1}{q}\right)\nm{u}^2\,-\,\mu\left(\frac{1}{p}-\frac{1}{q}\right)\nm{u}^p_p\,\geq\,0.
\end{equation}
When $\mu\leq0$, then $u=0$ follows immediately.
Accordingly, we need assume $\mu>0$ in general.
In this case, we correspondingly have
\begin{equation}
\mu\,c_2^p\left(\frac{1}{q}-\frac{1}{p}\right)\nm{u}^p\,\geq\,\mu\left(\frac{1}{q}-\frac{1}{p}\right)\nm{u}^p_p\,\geq\,\left(\frac{1}{q}-\frac{1}{2}\right)\nm{u}^2,\nonumber
\end{equation}
from which one can deduce that
\begin{equation}\label{e3.11}
\nm{u}\,\geq\,\left\{\mu^{-1}\,\frac{\left(\frac{1}{q}-\frac{1}{2}\right)}{c_2^p\left(\frac{1}{q}-\frac{1}{p}\right)}\right\}^{\frac{1}{p-2}}\,\to\,\infty
\end{equation}
when $\mu\to0^+$.
This finishes the proof for part {\bf(a)} of proposition \ref{P3.1}.

In addition, let $v\in\HarmU$ be such that $\varphi(v)\leq0$ and $\varphi^{\prime}(v)=0$.
Similarly, one has
\begin{equation}\label{e3.12}
\varphi(v)\,-\,\frac{1}{p}\,\varphi^{\prime}(v)(v)\,=\,\left(\frac{1}{2}-\frac{1}{p}\right)\nm{v}^2\,-\,\lambda\left(\frac{1}{q}-\frac{1}{p}\right)\nm{v}^q_q\,\leq\,0.
\end{equation}
When $\lambda\leq0$, then $v=0$ follows immediately.
Accordingly, we need assume $\lambda>0$ in general.
In this case, we correspondingly have
\begin{equation}
\lambda\,c_1^q\left(\frac{1}{q}-\frac{1}{p}\right)\nm{v}^q\,\geq\,\lambda\left(\frac{1}{q}-\frac{1}{p}\right)\nm{v}^q_q\,\geq\,\left(\frac{1}{2}-\frac{1}{p}\right)\nm{v}^2,\nonumber
\end{equation}
from which it yields readily that
\begin{equation}\label{e3.13}
\nm{v}\,\leq\,\left(\lambda\,\frac{c_1^q\left(\frac{1}{q}-\frac{1}{p}\right)}{\left(\frac{1}{2}-\frac{1}{p}\right)}\right)^{\frac{1}{2-q}}\,\to\,0
\end{equation}
when $\lambda\to0^+$.
This finishes the proof for part {\bf(b)} of proposition \ref{P3.1}.
\end{proof}

Finally, we consider the following analogous problem
\begin{equation}\label{e3.14}
\left\{\begin{array}{ll}
-\Delta\,u(x)\,+\,u(x)\,=\,0&\mathrm{in}\hspace{2mm}U,\\ \\
\Dnu(z)\,=\,\lambda\n{u(z)}^{q-2}u(z)\,+\,\mu\n{u(z)}^{p-2}u(z)&\mathrm{on}\hspace{2mm}\bdU,
\end{array}\right.
\end{equation}
in the standard Hilbert-Sobolev space $\HoU$, where all $u\in\HoU$ satisfy $u,\n{\gradu}\in\LtU$, and define the associated energy functional
\begin{equation}\label{e3.15}
\phi(u)\,:=\,\frac{1}{2}\InU\left[\,\n{\gradu}^2+\n{u}^2\,\right]\dx\,-\,\frac{\lambda}{q}\InbU\n{u}^q\dsg\,-\,\frac{\mu}{p}\InbU\n{u}^p\dsg.
\end{equation}
Note that in view of some result in \cite{AH1}, we have the following direct sum
\begin{equation}\label{e3.16}
\HoU\,=\,\HoUo\,\oplus\,\NrhU,
\end{equation}
where $\NrhU$ is the Hilbert subspace of $\HoU$ of all functions $u$ satisfying
\begin{equation}
\aip{u,v}_{1,\hspace{0.2mm}2}\,=\,\InU\left[\,\gradu\cdot\gradv+u\,v\,\right]\dx\,=\,0\hspace{7mm}\mathrm{for}\hspace{2mm}\mathrm{all}\hspace{2mm}v\in\CcOU\nonumber
\end{equation}
and $\HoUo$ is the closure of $\CcOU$ with respect to the standard $H^1$-norm.
Applying a similar procedure as in the proof of theorem \ref{T1.1}, we can obtain the following result.

\begin{thm}\label{T3.2}
Assume condition {\bf(B1)} holds and $1<q<2<p<2^*$.\\
{\bf(a)} When $\lambda\in\R$ and $\mu>0$ arbitrary, then problem \eqref{e3.14} has a sequence $\left\{\uk\right\}$ of solutions in $\NrhU$ such that $\phi(\uk)>0$ and $\phi(\uk)\to\infty$ as $k\to\infty$.\\
{\bf(b)} When $\lambda>0$ and $\mu\in\R$ arbitrary, then problem \eqref{e3.14} has a sequence $\left\{\vk\right\}$ of solutions in $\NrhU$ such that $\phi(\vk)<0$ and $\phi(\vk)\to0^-$ as $k\to\infty$.
\end{thm}

\vspace{2ex}
\section{$p$-Laplacian Steklov eigenvalue problems}\label{Steklov}
As mentioned earlier, the beauty of the paper \cite{AH2} is the discovery of the generalizations to high dimensions of the classical 3d \textsl{Laplace}'s \textsl{spherical harmonics} exterior to the unit ball: the \textsl{exterior harmonic Steklov eigenvalue problems} whose full spectra are derived there.
This section is devoted to the description of the exterior $p$-harmonic Steklove eigenvalue problems in the function space $\EpU$ when $N\geq3$ and $1<p<N$.
Similar results on bounded regions may be found in the interesting paper of Torn\'{e} \cite{To}.

Recall the \textsl{$p$-Laplacian} is defined as $\Delta_p\hspace{0.2mm}u:=\mathrm{div}\left(\n{\gradu}^{p-2}\gradu\right)$.
The \textsl{exterior $p$-harmonic Steklove eigenvalue problems} are to seek weak solutions of the problem below
\begin{equation}\label{e4.1}
\begin{split}
-\Delta_p\,u(x)&\,=\,0\hspace{2mm}\mathrm{in}\hspace{2mm}U,\\
\mathrm{subject}\hspace{2mm}\mathrm{to}\hspace{2mm}\n{\gradu(z)}^{p-2}\Dnu\,u(z)&\,=\,\delta\n{u(z)}^{p-2}u(z)\hspace{2mm}\mathrm{on}\hspace{2mm}\bdU,
\end{split}
\end{equation}
in $\EpU$.
This problem has been well-developed on bounded regions over a century since Stekloff \cite{St1, St2}, yet only been investigated recently in \cite{AH2, AH3, AH4, Ha1, Ha2, Ha4}.

Take $\nm{\cdot}:=\nm{\cdot}_{p,\hspace{0.2mm}\nabla}$ in this section.
Define two functionals on $\EpU$ as
\begin{equation}\label{e4.2}
\varphi(u)\,:=\,\frac{1}{p}\InbU\n{u}^p\dsg\hspace{4mm}\mathrm{and}\hspace{4mm}\psi(u)\,:=\,\frac{1}{p}\InU\n{\gradu}^p\dx.
\end{equation}
Accordingly, given $u\in\EpU$, write two linear functionals on $\EpU$ to be
\begin{equation}\label{e4.3}
\left\{\begin{array}{ll}
\Pc_u(v)\,:=\,p\,\varphi(u)\,\psi^{\prime}(u)(v)\\ \\
\Bc_u(v)\,:=\,\varphi^{\prime}(u)(v)\,-\,\Pc_u(v)
\end{array}\right.\hspace{7mm}\mathrm{for}\hspace{2mm}\mathrm{each}\hspace{2mm}v\in\EpU.
\end{equation}

Since $\EpU$ is a reflexive, uniformly convex Banach space when $p>1$ (see \cite{AH3}), there exists a unique element, say, $u_{\Bc}$ in $\EpU$, from \textsl{Riesz}'s \textsl{theorem}, such that
\begin{equation}\label{e4.4}
\Bc_u(u_{\Bc})\,=\,\nm{\Bc_u}^2_*\,=\,\nm{u_{\Bc}}^2;
\end{equation}
therefore, one finds a homeomorphism $\mathfrak{H}:\EpU\to\EpU$ such as $\mathfrak{H}(u):=u_{\Bc}$.
Noticing $\mathfrak{H}$ is odd, and bounded, uniformly continuous on the set $\mathbf{S_1}:=\left\{u\in\EpU:\nm{u}=1\right\}$, there are constants $t_0,\upsilon_1,\upsilon_0>0$ such that, for all $t\in\left[-t_0,t_0\right]$ and $u\in\mathbf{S}_1$, one has
\begin{equation}\label{e4.5}
\upsilon_1\,\geq\,\nm{u+tu_{\Bc}}\,\geq\,\upsilon_0.
\end{equation}

Now, define a flow $H:\mathbf{S}_1\times\left[-t_0,t_0\right]\to\mathbf{S}_1$ by
\begin{equation}\label{e4.6}
H(u,t)\,:=\,\frac{u+tu_{\Bc}}{\nm{u+tu_{\Bc}}}.
\end{equation}
Then, $H$ is odd in $u$, with $H(u,0)=u$, uniformly continuous and verifies the property below.

\begin{lem}\label{L4.1}
There exists a map $\ell(u,t):\mathbf{S}_1\times\left(-t_0,t_0\right)\to\R$ such that $\ell(u,t)\to0$ as $t\to0$, uniformly on $\mathbf{S}_1$, and, for each $u\in\mathbf{S}_1$ and $t\in\left(-t_0,t_0\right)$, we have
\begin{equation}\label{e4.7}
\varphi(H(u,t))\,=\,\varphi(u)\,+\,\int_0^t\left[\,\nm{u_{\Bc}}^2+\ell(u,s)\,\right]\ds.
\end{equation}
\end{lem}

\begin{proof}
Since
\begin{equation}\label{e4.8}
\varphi(H(u,t))\,=\,\varphi(u)\,+\,\int_0^t\,\varphi^{\prime}(H(u,s))\left(\frac{\pal H(u,s)}{\pal s}\right)\ds,
\end{equation}
we can define, in view of \eqref{e4.4},
\begin{equation}\label{e4.9}
\ell(u,t)\,:=\,\varphi^{\prime}(H(u,t))\left(\frac{\pal H(u,t)}{\pal t}\right)\,-\,\Bc_{u}(u_{\Bc}),
\end{equation}
and have $\ell(u,0)=0$ uniformly on $\mathbf{S}_1$ via \eqref{e4.5}.
Actually, a routine computation leads to
\begin{equation}
\frac{\pal H(u,t)}{\pal t}\,=\,\frac{u_{\Bc}}{\nm{u+tu_{\Bc}}}\,-\,\frac{u+tu_{\Bc}}{\nm{u+tu_{\Bc}}^{p+1}}\InU\n{\nabla(u+tu_{\Bc})}^{p-2}\nabla(u+tu_{\Bc})\cdot\gradu_{\Bc}\,\dx,\nonumber
\end{equation}
from which one deduces easily that, noticing that $\nm{u}=1$ on $\mathbf{S}_1$,
\begin{equation}\label{e4.10}
\frac{\pal H(u,0)}{\pal t}\,=\,u_{\Bc}\,-\,u\InU\n{\gradu}^{p-2}\gradu\cdot\gradu_{\Bc}\,\dx;
\end{equation}
this further implies that, remembering the fact $H(u,0)=u$ now,
\begin{equation}
\begin{split}
\varphi^{\prime}(H(u,0))\left(\frac{\pal H(u,0)}{\pal t}\right)&\,=\,\varphi^{\prime}(u)(u_{\Bc})\,-\,\varphi^{\prime}(u)(u)\,\psi^{\prime}(u)(u_{\Bc})\,=\,\varphi^{\prime}(u)(u_{\Bc})\,-\,\Pc_u(u_{\Bc}),\nonumber
\end{split}
\end{equation}
which together with \eqref{e4.3} and \eqref{e4.9} gives the desired result as $\ell$ is bounded by \eqref{e4.5}.
\end{proof}

Using this result, we can derive a version of \textsl{deformation lemma}.

\begin{prop} \label{P4.2}
Given a constant $\kappa>0$, suppose there are constants $\varsigma>0$ and $\tau\in\left(0,\kappa\right)$, such that $\nm{u_{\Bc}}\geq\varsigma$ on $\mathtt{V}_{\tau}:=\left\{u\in\mathbf{S}_1:\n{\varphi(u)-\kappa}\leq\tau\right\}$.
Then, for every compact, symmetric subset $\mathbf{G}$ of $\mathbf{S}_1$, one finds a constant $\epsilon\in\left(0,\tau\right)$ and an associated odd map $H_{\epsilon}:\mathbf{S}_1\to\mathbf{S}_1$ that is continuous on $\mathtt{V}_{\epsilon}\cap\mathbf{G}$ and $H_{\epsilon}(\mathtt{V}_{\epsilon}\cap\mathbf{G})\subseteq\varphi_{\kappa+\epsilon}$, where $\varphi_{\kappa+\epsilon}:=\left\{u\in\mathbf{S}_1:\varphi(u)\geq\kappa+\epsilon\right\}$.
\end{prop}

\begin{proof}
As $\ell(u,0)=0$ uniformly on $\mathbf{S}_1$, we can choose $t_1\in\left(0,t_0\right)$ such that $\n{\ell(u,t)}\leq\varsigma^2/2$ for each $u\in\mathbf{G}\varsubsetneq\mathbf{S}_1$ and $t\in\left[-t_1,t_1\right]$ since $\mathbf{G}$ is compact.
Write $\epsilon:=\min\left\{\tau,\varsigma^2t_1/4\right\}$.
Then, \eqref{e4.4} and \eqref{e4.7} implies that, for every $u\in\mathtt{V}_{\epsilon}\cap\mathbf{G}\subseteq\mathtt{V}_{\tau}$, $\varphi(H(u,t_1))\geq\kappa-\epsilon+\frac{\varsigma^2t_1}{2}\geq\kappa+\epsilon$; that is, $H(u,t_1)\in\varphi_{\kappa+\epsilon}$ for all $u\in\mathtt{V}_{\epsilon}\cap\mathbf{G}$.
As such, we define an odd map via \eqref{e4.6} that is continuous on $\mathtt{V}_{\epsilon}\cap\mathbf{G}$ by (because $\varphi$ is even and $\mathbf{G}=-\,\mathbf{G}$ so that $\mathtt{V}_{\epsilon}\cap\mathbf{G}$ again is symmetric)
\begin{equation}\label{e4.11}
H_{\epsilon}(u)\,:=\,\left\{\begin{array}{ll}
H(u,t_1)&\mathrm{when}\hspace{2mm}u\in\mathtt{V}_{\epsilon}\cap\mathbf{G},\\ \\
u&\mathrm{when}\hspace{2mm}u\in\mathbf{S}_1\setminus\mathtt{V}_{\epsilon}\cap\mathbf{G},
\end{array}\right.
\end{equation}
and have $H_{\epsilon}(\mathtt{V}_{\epsilon}\cap\mathbf{G})\subseteq\varphi_{\kappa+\epsilon}$, as claimed.
This in turn completes the proof.
\end{proof}

Next, define
\begin{equation}\label{e4.12}
\kappa_n\,:=\,\sup\limits_{\mathbf{G}\in\mathscr{G}_n}\,\min\limits_{u\in\mathbf{G}}\,\varphi(u)\,\geq\,0.
\end{equation}
Here, $\mathscr{G}_n:=\left\{\mathbf{G}\varsubsetneq\mathbf{S}_1:\mathbf{G}\hspace{1.2mm}\mathrm{compact},\hspace{1.2mm}\mathbf{G}=-\mathbf{G}\hspace{1.2mm}\mathrm{and}\hspace{1.2mm}\gamma(\mathbf{G})\leq n\right\}$, with $\gamma$ being the \textsl{Krasnoselskii genus} (see, for example, section ii.5 of \cite{St} for more detailed descriptions).

Then, we can prove the following main results.

\begin{thm}\label{T4.3}
For every $n\geq1$, $\kappa_n>0$ and there exists a function $s_n\in\EpU$ such that $\varphi(s_n)=\kappa_n$; in addition, $s_n$ is a weak solution of \eqref{e4.1} with $\delta=\delta_n:=\frac{1}{p\,\kappa_n}>0$.
\end{thm}

\begin{proof}
As $\gamma(\mathbf{S}_1)=\infty$, $\kappa_n$ is well-defined in the sense that $\mathscr{G}_n\neq\emptyset$ for each $n\in\mathbb{N}$.
Select a set $\mathbf{G}_n\in\mathscr{G}_n$ with $u\neq0$ $\sigma$ \textit{a.e.} on $\bdU$ for all $u\in\mathbf{G}_n$ to derive $\kappa_n\geq\min\limits_{u\in\mathbf{G}_n}\varphi(u)>0$.

Next, given $n\geq1$, there exists a sequence $\left\{u_{n,\hspace{0.2mm}k}\right\}$ in $\mathbf{S}_1$ such that $\varphi(u_{n,\hspace{0.2mm}k})\to\kappa_n$.
Using a subsequence if necessary, it implies $u_{n,\hspace{0.2mm}k}\rightharpoonup s_n\in\EpU$ yet $u_{n,\hspace{0.2mm}k}\to s_n\in\LpbU$ in view of result 2.1, when $k\to\infty$.
Thus, one deduces that $\varphi(u_{n,\hspace{0.2mm}k})\to\varphi(s_n)=\kappa_n$.

Moreover, as linear functionals on $\EpU$, $\Bc_{u_{n,\hspace{0.2mm}k}}\to0$ when $k\to\infty$.
First, by definition of $\kappa_n$, one can find a set $\tilde{\mathbf{G}}_n\in\mathscr{G}_n$, with $\kappa_n-\epsilon\leq\varphi(u)\leq\kappa_n+\epsilon$, for each $u\in\tilde{\mathbf{G}}_n$ and some suitably small $\epsilon\in\left(0,\frac{\kappa_n}{4}\right)$; now, if we suppose on the contrary $\nm{\Bc_{u}}_*>\varsigma>0$ uniformly on $\left\{u\in\mathbf{S}_1:\frac{\kappa_n}{2}\leq\varphi(u)\leq\frac{3\kappa_n}{2}\right\}$, proposition \ref{P4.2} provides us with a continuous, odd map $H_{\epsilon}$ on $\tilde{\mathbf{G}}_n$ such that $H_{\epsilon}\big(\tilde{\mathbf{G}}_n\big)\in\mathscr{G}_n$ and $H_{\epsilon}\big(\tilde{\mathbf{G}}_n\big)\subseteq\varphi_{\kappa_n+\epsilon}$ - a contradiction thus arrives.

Since $\Bc_{u_{n,\hspace{0.2mm}k}}\to0$ when $k\to\infty$, one infers that, in view of \eqref{e4.2} and \eqref{e4.3},
\begin{equation}
\lim_{k\to\infty}\InU\n{\gradu_{n,\hspace{0.2mm}k}}^{p-2}\gradu_{n,\hspace{0.2mm}k}\cdot\gradv\,\dx\,=\,\Fc(v)\hspace{7mm}\mathrm{for}\hspace{2mm}\mathrm{all}\hspace{2mm}v\in\EpU,\nonumber
\end{equation}
where $\Fc(v):=\frac{1}{p\,\kappa_n}\InbU\n{s_n}^{p-2}s_nv\,\dsg$ is a linear functional on $\EpU$.
Using result 2.1 again, plus $u_{n,\hspace{0.2mm}k}\rightharpoonup s_n\in\EpU$ yet $u_{n,\hspace{0.2mm}k}\to s_n\in\LpbU$ as $k\to\infty$, we obtain
\begin{equation}\label{e4.13}
\InU\n{\grads_n}^{p-2}\grads_n\cdot\gradv\,\dx\,=\,\delta_n\InbU\n{s_n}^{p-2}s_n\,v\,\dsg\hspace{7mm}\mathrm{for}\hspace{2mm}\mathrm{all}\hspace{2mm}v\in\EpU.
\end{equation}
Here, $\delta_n:=\frac{1}{p\,\kappa_n}>0$.
As such, $s_n\in\EpU$ is a weak solution of problem \eqref{e4.1}.
\end{proof}

\begin{thm}\label{T4.4}
Define $\kappa_n$ as in \eqref{e4.12} and $\delta_n$ by $\frac{1}{p\,\kappa_n}$ for each $n\in\mathbb{N}$.
Then, one has
\begin{equation}\label{e4.14}
\lim_{n\to\infty}\,\delta_n\,=\,\infty.
\end{equation}
\end{thm}

\begin{proof}
The conclusion \eqref{e4.14} follows if we can show that $\lim\limits_{n\to\infty}\kappa_n=0$.

For each $n\geq1$, choose $\mathcal{E}_n$ to be a linear subspace of $\EpU$ of dimension $n$ such that $u\neq0$ $\sigma$ \textit{a.e.} on $\bdU$ for every $u\in\mathcal{E}_n$, and denote its complement in $\EpU$ by $\mathcal{E}^{c}_n$.
Without loss of generality, assume further that $\mathcal{E}_1\varsubsetneq\mathcal{E}_2\varsubsetneq\cdots\varsubsetneq\mathcal{E}_n\varsubsetneq\cdots\varsubsetneq\EpU$.
Note that \eqref{e4.13} guarantees our choice as clearly $s_n|_{\bdU}\neq0$ $\sigma$ \textit{a.e.} on $\bdU$ for all $n\in\mathbb{N}$.
Also, we have
\begin{equation}\label{e4.15}
\overline{\cup_{n\in\mathbb{N}}\,\mathcal{E}_n}\,\dot\bigcup\,\EopU\,=\,\EpU,
\end{equation}
where $\EopU$ denotes the subspace of $\EpU$ that is the closure of $\CcOU$ with respect to the gradient $L^p$-norm \eqref{e2.1} and the notation $\dot\cup$ means disjoint union.

Now, define $\tilde{\kappa}_n:=\sup\limits_{\mathbf{G}\in\mathscr{G}_n}\min\limits_{u\in\mathbf{G}\cap\mathcal{E}^{c}_n}\varphi(u)$ to give $\tilde{\kappa}_n\geq\kappa_n>0$.
Then, one proves $\lim\limits_{n\to\infty}\tilde{\kappa}_n=0$.
Actually, if not, there is a constant $\varepsilon>0$ such that $\tilde{\kappa}_n\geq\varepsilon$ for all $n\geq1$.
Thus, a set $\breve{\mathbf{G}}_n\in\mathscr{G}_n$ exists such that $\tilde{\kappa}_n\geq\min\limits_{u\in\breve{\mathbf{G}}_n\cap\mathcal{E}^{c}_n}\varphi(u)\geq\frac{\varepsilon}{2}>0$ for each $n\in\mathbb{N}$, so that we find a sequence $\left\{u_n\right\}$, with $u_n\in\breve{\mathbf{G}}_n\cap\mathcal{E}^{c}_n$, satisfying $\varphi(u_n)\geq\frac{\varepsilon}{2}$ uniformly.
Keep in mind $\breve{\mathbf{G}}_n\varsubsetneq\mathbf{S}_1$; from \eqref{e4.15} and resorting to a subsequence if necessary, one has $u_n\to0\in\LpbU$, and thus $\varphi(u_n)\to0$, as $n\to\infty$.
A contradiction follows and thereby one finishes the proof.
\end{proof}

Finally, it is worth to mention here that the following problem
\begin{equation}\label{e4.16}
\begin{split}
-\Delta_p\,u(x)\,+\,\n{u(x)}^{p-2}u(x)&\,=\,0\hspace{2mm}\mathrm{in}\hspace{2mm}U,\\
\mathrm{subject}\hspace{2mm}\mathrm{to}\hspace{2mm}\n{\gradu(z)}^{p-2}\Dnu\,u(z)&\,=\,\delta\n{u(z)}^{p-2}u(z)\hspace{2mm}\mathrm{on}\hspace{2mm}\bdU,
\end{split}
\end{equation}
can be studied analogically in the space $\WpU$ by use of some results in \cite{AH1}.


\end{document}